\documentclass[a4paper,11pt]{article}

\usepackage{graphicx}
\usepackage{amsmath}
\usepackage{amsfonts}
\usepackage{titlesec}
\usepackage{amsthm}
\usepackage{amssymb}
\usepackage{color}
\usepackage{mathtools}
\usepackage{geometry}
\usepackage{setspace}
\usepackage{times}
\usepackage{algorithm2e}

\usepackage{verbatim}

\numberwithin{equation}{section}

\newtheorem{theorem}{Theorem}[section]

\theoremstyle{definition}

\newtheorem{remark}[theorem]{Remark}

\definecolor{light-gray}{gray}{0.69}
\definecolor{light-red}{rgb}{1.0,0.4,0.4}

\definecolor{light-blue}{rgb}{0.4,0.45,1}
\definecolor{light-green}{rgb}{0.5,0.8,0.0}
\definecolor{dark-green}{rgb}{0.0,0.4,0.0}
\definecolor{dark-red}{rgb}{1.0,0.3,0.3}

\definecolor{dark-gray}{gray}{0.59}
\definecolor{very-dark-gray}{gray}{0.39}
\definecolor{lighter-red}{rgb}{1.0,0.6,0.6}

\definecolor{ocker_hell}{rgb}{0.75,0.7,0.4}
\definecolor{gelb_dunkel}{rgb}{0.75,0.7,0.0}
\definecolor{gruen_hell}{rgb}{0.5,0.8,0.0}

\definecolor{dark-blue}{rgb}{0.0,0.0,0.5}
\definecolor{new-blue}{rgb}{0.0,0.0,0.8}

\definecolor{lila}{rgb}{0.5,0.0,0.5}
\definecolor{dark-red}{rgb}{0.5,0.0,0.0}

\begin{document}

\begin{center}
{\huge Dispersion-minimizing quadrature rules for $C^1$ quadratic isogeometric analysis }
\end{center}

\renewcommand{\thefootnote}{\fnsymbol{footnote}}
\renewcommand{\thefootnote}{\arabic{footnote}}

\begin{center}
{\large Quanling Deng$^{a,}$\footnote{Corresponding author. qdeng12@gmail.com}, Michael Barto\v{n}$^b$, Vladimir Puzyrev$^a$,  Victor Calo$^{a,c}$} \\
$^a${\small Department of Applied Geology, Curtin University, Kent Street, Bentley, Perth, WA 6102, Australia} \\
$^b${\small Basque Center for Applied Mathematics, Alameda de Mazarredo 14, 48009 Bilbao, Basque Country, Spain} \\
$^c${\small Mineral Resources, CSIRO, Kensington, Perth, WA 6152, Australia}
\end{center}

\begin{abstract}
We develop quadrature rules for the isogeometric analysis of wave propagation and structural vibrations that minimize the discrete dispersion error of the approximation. The rules are optimal in the sense that they only require two quadrature points per element to minimize the dispersion error  \cite{bartovn2017generalization}, and they are equivalent to the optimized blending rules we recently described. Our approach further simplifies the numerical integration:  instead of blending two three-point standard quadrature rules, we construct directly a single two-point quadrature rule that reduces the dispersion error  to the same order for uniform meshes with periodic boundary conditions. Also, we present a 2.5-point rule for both uniform and non-uniform meshes with arbitrary boundary conditions. Consequently, we reduce the computational cost by using the proposed quadrature rules. Various numerical examples demonstrate the performance of these quadrature rules.
\end{abstract}


\section{Introduction} \label{sec:intr} 

Quadrature rules play an important role in the implementation of various numerical methods for solving partial differential equations. Fewer quadrature points result in a lower computational cost, however, the reduction of the quadrature points should not reduce the quality of the approximation. The design of efficient quadrature rules for isogeometric analysis (see Hughes \textit{et al.} \cite{hughes2005isogeometric, bazilevs2006isogeometric, cottrell2006isogeometric, cottrell2009isogeometric}) is of interest as the continuity properties of the spline basis functions may require fewer quadrature points. The quadrature rules should preserve the optimal convergence of the numerical approximation to the exact solution. Traditionally, Gauss rules for discontinuous polynomial spaces are used, however, these choices are far from being optimal in general \cite{auricchio2012simple}.

The construction of efficient quadrature rules for isogeometric analysis was initially considered by Hughes \textit{et al.} in \cite{hughes2010efficient} in 2010. Taking advantage of the smoothness of the basis functions across element boundaries, a half-point rule that is independent of the polynomial order of the basis functions was developed.  The new rule has advantages when compared to the traditional ones. The rule is optimal as it exactly integrates the spline basis functions with the minimum of number of quadrature sampling points. The rule is designed for uniform univariate splines and is Gaussian, that is, optimal in the sense of the minimum number of quadrature points. However, the rule is exact only for infinite domains or for the spline spaces that have a special structure at the boundaries of finite domains.  To make the rule exact for a general spline space over finite domains, additional quadrature points are introduced at the boundary elements, resulting in nearly-optimal quadrature rules [5]. These  non-Gaussian rules come as solutions of non-linear, possibly ill-conditioned, systems and possess both positive and negative weights.

Other works in this direction are reported in \cite{schillinger2014reduced,hiemstra2016optimal,bartovn2016gaussian,bartovn2016optimal,bartovn2017gauss}.  Optimal and reduced quadrature rules for tensor product and hierarchically refined splines for isogeometric analysis were developed in \cite{hiemstra2016optimal}. Gaussian rules for spline spaces of various degrees and continuities were derived in \cite{bartovn2016optimal,bartovn2017gauss}. Using the homotopy continuation argument [9], Gaussian rules can be derived by continuously modifying the spline space (knot vector) and by tracing numerically the rule, which is given by solving a certain algebraic system. These rules guarantee exactness of the integration up to machine precision, and the property of being Gaussian also directly implies that all weights are positive \cite{micchelli1977moment}.

The study of dispersion error minimization for isogeometric analysis was initially studied numerically in Puzyrev \textit{et al.} \cite{puzyrev2017dispersion} and analytically in Calo \textit{et al.} \cite{calo2017dispersion}. For general dispersion analysis of isogeometric discretizations, we refer the readers to \cite{hughes2008duality,hughes2014finite} and the references therein. Particularly, in Hughes \textit{et al.} \cite{hughes2008duality}, a duality principle between the dispersion analysis and the spectral analysis was established and the analysis unified.  

The study of dispersion analysis of the finite element method has a rich literature; see for example Thomson and Pinsky\cite{thompson1994complex,thompson1995galerkin}, Ihlenburg and Babuska \cite{ihlenburg1995dispersion}, Ainsworth \cite{ainsworth2004discrete, ainsworth2009dispersive,ainsworth2010optimally}, and others \cite{harari1997reducing,harari2000analytical,he2011dispersion}. Thomson and Pinsky studied the dispersive effects of the finite element methods with Legendre, spectral, and Fourier local approximation basis for the Helmholtz equation in \cite{thompson1994complex}. They found that the choice of basis functions had a negligible effect on the dispersion errors. This is due to the low continuity ($C^0$ continuity) of the basis functions. Hughes \textit{et al.} \cite{hughes2008duality} showed that the dispersion error of the isogeometric analysis with high continuity (up to $C^{p-1}$ for $p$-th order basis function) on the basis functions is smaller than that of the lower continuity finite element counterparts.

The $2p$-optimal convergence rate of the dispersion error for the $p$-th order standard finite elements was established in \cite{ainsworth2004discrete}. In 2009, Ainsworth and Wajid  \cite{ainsworth2009dispersive} extended this analysis to arbitrary spectral element methods. Based on Marfurt's conjecture \cite{marfurt1984accuracy} that the most promising and efficient method for computing wave propagation is to blend the finite element method with the spectral element method with appropriate weights,  Ainsworth and Wajid beautifully established the optimal blending of these two methods in \cite{ainsworth2010optimally}. A superconvergence (order $2p+2$ for $p$-th order polynomial approximation) result was obtained for arbitrary order of polynomial approximation, which includes the fourth order superconvergence result obtained by a modified integration rule for linear finite elements in \cite{guddati2004modified}.

To the best of our knowledge, this is the first paper studying the design of optimal quadrature rules which  minimize the dispersion errors of the isogeometric analysis for the wave propagation and structural vibration problems. The dispersion error-minimizing quadratures, that combine Gauss-Legendre and Gauss-Lobatto rules proposed in \cite{puzyrev2017dispersion,calo2017dispersion}, are not  efficient as the two traditional quadrature rules are used for each integration evaluation. Herein, we design quadrature rules that minimize the dispersion error and minimize the number of quadrature points. A rule that has minimal number of evaluation points per element (two in the case of a uniform mesh with periodic boundary conditions) is the solution of a non-linear system of algebraic equations which, due to the low polynomial degree, admits a closed form formula.  We also design a quadrature rule that minimizes dispersion for the larger $C^0$ quadratic space, which leads to a quadrature rule that uses 2.5 points per element as it exactly integrates discontinuous cubic functions on  the mesh. This rule is effective for non-uniform meshes and arbitrary boundary conditions.

The rest of this paper is organized as follows. Section \ref{sec:ps} describes the isogeometric discretization of an eigenvalue problem. In Section \ref{sec:ddr}, we present the constraints minimizing both dispersion error and the number of quadrature points and set up the equations for the quadrature weights and points. Both two-point and 2.5-point rules are considered here. Section \ref{sec:num} studies numerical examples to demonstrate the performance of both the two-point and 2.5-point rules. Concluding remarks are given in Section \ref{sec:conclusion}.

\section{Model problem and its discretization} \label{sec:ps}
This section follows closely \cite{calo2017dispersion}. In order to illustrate the main ideas, we consider stationary waves as described by the Helmholtz equation
\begin{equation} \label{eq:pde}
\begin{aligned}
\Delta u +  \omega^2 u & = 0 \quad  \text{in} \quad \Omega, \\
u & = 0 \quad \text{on} \quad \partial \Omega,
\end{aligned}
\end{equation}
where $\Omega = [0, 1]^d \subset \mathbb{R}^d, d= 1,2,3, \Delta = \nabla^2$ is the Laplacian and $\omega = \omega_f/c$ with $\omega_f$ being the frequency of a particular sinusoidal wave and $c$ being the speed of sound of the medium.

For an open set $S \subset \mathbb{R}^d $  with Lipschitz boundary, we denote by $W^{m,p}(S)$ a Sobolev space equipped with the norm $\| \cdot \|_{m,p,S}$ and the semi-norm $| \cdot |_{m,p,S}$ where $m$ is the weak derivative order and $p$ corresponds to the $p$ in $L^p$ space. We use standard notation. If $p=2$, we omit $p$ and utilize $H^m(S)$  for Hilbert spaces and $H^m_0(S)$ for Hilbert spaces with functions vanishing at the boundary for $m>0$. 
The variational formulation of \eqref{eq:pde} is to find $u \in H^1_0(\Omega)$ such that 
\begin{equation} \label{eq:vf}
B(u, v) = 0 \quad \forall \ v \in H^1_0(\Omega), 
\end{equation}
where 
\begin{equation}
B(w, v) = a(w, v) - \omega^2 b(w, v)
\end{equation}
with
$
a(w, v) = (\nabla w, \nabla v)
$ and $
b(w, v) = (w, v)
$. Here $(\cdot, \cdot)$ is the $L^2$ inner product.

Let $\mathcal{T}_h$ be a discretization of the bounded open domain $\Omega$ and we denote each element as $K$ such that $\bar\Omega = \cup_{K\in \mathcal{T}_h}  K$. Let $h = \max_{K\in \mathcal{T}_h} \text{diameter}(K)$.
The Galerkin-type numerical methods seek $u_h \in V_h$ such that 
\begin{equation} \label{eq:vfh}
B(u_h, v_h) =  0 \quad \forall \ v_h \in V_h. 
\end{equation}

Different trial spaces $V_h$ lead to different numerical methods. We focus on isogeometric analysis in this work. We denote $\phi_a = \phi_a(x)$ the  B-spline basis functions  we use in isogeometric analysis. Then $ V_h = \text{span} \{\phi_a\}$.  

In practice, the integrals involved in $a(u_h, v_h) $ and $b(u_h, v_h)$ are evaluated numerically, that is, approximated by quadrature rules. On a reference element $\hat K$, a quadrature rule is of the form
\begin{equation} \label{eq:qr}
\int_{\hat K} \hat f(\hat{\boldsymbol{x}}) \ \text{d} \hat{\boldsymbol{x}} \approx \sum_{l=1}^{N_q} \hat{\varpi}_l \hat f (\hat{n_l}),
\end{equation}
where $\hat{\varpi}_l$ are the weights, $\hat{n_l}$ are the nodes, and $N_q$ is the number of nodes. For each element $K$, we assume that there is an invertible affine map $\sigma$ such that $K = \sigma(\hat K)$, which leads to the correspondence between the functions on $K$ and $\hat K$. Assuming $J_K$ is the corresponding Jacobian of the mapping, \eqref{eq:qr} induces a quadrature rule over the element $K$ given by
\begin{equation} \label{eq:q}
\int_{K}  f(\boldsymbol{x}) \ \text{d} \boldsymbol{x} \approx \sum_{l=1}^{N_q} \varpi_{l,K} f (n_{l,K}),
\end{equation}
where $\varpi_{l,K} = \text{det}(J_K) \hat \varpi_l$ and $n_{l,K} = \sigma(\hat n_l)$.

Applying quadrature rules to \eqref{eq:vfh}, we have the approximate form
\begin{equation} \label{eq:vfho}
\tilde B_h(u_h, v_h) = \tilde a_h(\tilde u_h, v_h) -  \omega^2 \tilde b_h(\tilde u_h, v_h) \quad \forall \ v_h \in V_h,
\end{equation}
where 
\begin{equation} \label{eq:ba}
\tilde a_h(w, v) = \sum_{K \in \mathcal{T}_h} \sum_{l=1}^{N_q} \varpi_{l,K}^{(1)} \nabla w (n_{l,K}^{(1)} ) \cdot \nabla v (n_{l,K}^{(1)} ),
\end{equation}
and
\begin{equation} \label{eq:bb}
\tilde b_h(w, v) = \sum_{K \in \mathcal{T}_h} \sum_{l=1}^{N_q} \varpi_{l,K}^{(2)} w (n_{l,K}^{(2)} ) v (n_{l,K}^{(2)} ),
\end{equation}
where $\{\varpi_{l,K}^{(1)}, n_{l,K}^{(1)} \}$ and $\{\varpi_{l,K}^{(2)}, n_{l,K}^{(2)} \}$ specify two (possibly different) quadrature rules.

With quadrature rules, if we substitute the basis functions of $V_h$ into \eqref{eq:vfho}, this leads to the linear algebra problem
\begin{equation} \label{eq:las}
( \mathbf{K} -\omega^2 \mathbf{M} ) \mathbf{U} = \mathbf{0}
\end{equation}
where $\mathbf{K}$ and $\mathbf{M}$ are the global stiffness and mass matrices with entries $\mathbf{K}_{ab} = \tilde a_h(\phi_a, \phi_b)$, $\mathbf{M}_{ab} = \tilde b_h(\phi_a, \phi_b),$ and $\mathbf{U}$ is the unknown vector. This system \eqref{eq:las} admits a nontrivial solution in the view of \cite{ainsworth2004discrete}.

\section{Quadratures that minimize the dispersion error }\label{sec:ddr}

For simplicity, we consider a one-dimensional problem where $\Omega = \mathbb{R}$. We describe the general framework for all orders of isogeometric analysis and then focus on the $C^1$ quadratic case as follows. 

We denote by $C^k_p$ the space of piecewise polynomials of order $p$ and continuity $k$. Let $\phi_a = \phi_a(x)$ be a B-spline basis function of isogeometric analysis of order $p$ with maximum continuity $C^{p-1}$, thus $\phi_a \in C^{p-1}_p$. The functions for the integration corresponding to the stiffness $\tilde a_h(\phi_a, \phi_b)$ are in the space $C^{p-2}_{2p-2}$ while the functions for the integration corresponding to the mass $\tilde b_h(\phi_a, \phi_b)$ are in the space $C^{p-1}_{2p}$. Thus, to integrate both stiffness and mass matrices exactly, one needs quadrature rules which integrate all the functions in the space $C^{p-2}_{2p}$ exactly.

An $m$ point Gauss-Legendre quadrature rule, denoted by $G_m$, integrates the space $C^{-1}_{2m-1}$ exactly, while an $m$ point Gauss-Lobatto quadrature rule, denoted by $L_m$, integrates the space $C^{-1}_{2m-3}$ exactly. For $p$-th order isogeometric elements, to integrate both the stiffness and mass matrices exactly, the rule $G_{p+1}$ is enough since $C^{p-2}_{2p}$ is a subset of $C^{-1}_{2p+1}$ but $G_p$ is not enough. Thus, an optimized rule (minimized number of quadrature points) can be developed by considering the $p-2$ order of continuity; see \cite{bartovn2016optimal}.

In the view of analysis in \cite{calo2017dispersion} in a one-dimensional setting, the stiffness matrix is integrated exactly if the space of $C_{2p-1}^{p-2}$  is fully integrated for $p$-th order isogeometric elements while the mass matrix can be under-integrated to minimize the dispersion errors. Our aim is to develop a quadrature rule that minimizes the dispersion errors.

Below, we focus on isogeometric analysis with $C^1$ quadratic B-spline basis functions.
Let $\phi_a = \phi_a(x)$ be a $C^1$ quadratic B-spline basis function. We seek an approximation of the form 
\begin{equation} \label{eq:1dsol}
U (x) = \sum_{a \in \mathbb{Z}} U^a \phi_a(x)
\end{equation}
satisfying 
\begin{equation} \label{eq:bho}
\tilde B_h(U, v_h) = 0.
\end{equation}

The quadrature rules we use to integrate \eqref{eq:ba} and \eqref{eq:bb} are the same. In the one-dimension case, we seek a quadrature rule which integrates the stiffness matrix exactly.  We denote our new quadrature rule $NQ_2$ for $C^1$ quadratic isogeometric elements with nodes $n_1, n_2$ in the reference interval $[0, 1]$ and weights $\varpi_1, \varpi_2$.

We apply this quadrature rule $NQ_2$ to \eqref{eq:bho} to obtain the following equation for the value $U^j$ of the approximation at node $x_j = jh, j \in \mathbb{Z}$
\begin{equation} \label{eq:nq2iga}
\begin{aligned}
\big(K_2 - \Lambda^2 M_2 \big) ( U^{j-2} + U^{j+2}) - \big(K_1 + \Lambda^2 M_1 \big) ( U^{j-1} + U^{j+1} ) + \big(K_0 - \Lambda^2 M_0 \big)  U^j = 0,
\end{aligned}  
\end{equation}
where $\Lambda = \omega h$ and
\begin{equation}
\begin{aligned}
K_0 & = 2( 3 n_1^2 \varpi_1-3 n_1 \varpi_1+3 n_2^2 \varpi_2-3 n_2 \varpi_2+\varpi_2+\varpi_1), \\
K_1 & = (1-2 n_1)^2 \varpi_1 + (1-2 n_2)^2 \varpi_2, \\
K_2 & = (n_1-1) n_1 \varpi_1+(n_2-1) n_2 \varpi_2, \\
M_0 & = \frac{1}{2} \Big( (3 n_1^4-6 n_1^3+3 n_1^2+1 ) \varpi_1+ (3 n_2^4-6 n_2^3+3 n_2^2+1) \varpi_2 \Big), \\
M_1 & = \frac{1}{4} \Big( (-4 n_1^4+8 n_1^3-4 n_1^2+1 ) \varpi_1+ (-4 n_2^4+8 n_2^3-4 n_2^2+1 ) \varpi_2 \Big), \\
M_2 & = \frac{1}{4} \Big( (n_1-1)^2 n_1^2 \varpi_1+(n_2-1)^2 n_2^2 \varpi_2 \Big).
\end{aligned}  
\end{equation}

We assume that the equation admits nontrivial Bloch wave \cite{odeh1964partial} solutions in the form $U^j = e^{ij \mu^{(2)}_{Q}h }$, where $i^2 = -1$ and the subindex $Q$ denotes the corresponding numerical quadrature, then \eqref{eq:nq2iga} simplifies to 
\begin{equation} \label{eq:nq2igadr}
\begin{aligned}
2 \big(K_2 - \Lambda^2 M_2 \big) \cos(2\mu^{(2)}_{NQ_2}h) - 2 \big(K_1 + \Lambda^2 M_1 \big) \cos(\mu^{(2)}_{NQ_2}h) + \big(K_0 - \Lambda^2 M_0 \big) = 0,
\end{aligned}  
\end{equation}
which is known as the \textit{discrete dispersion relation}  for the discrete method with a particular quadrature rule. Solving \eqref{eq:nq2igadr} for $\mu^{(2)}_{NQ_2}h$ and writing the expression as a series in $\Lambda$, we obtain
\begin{equation} \label{eq:nq2igade}
\begin{aligned}
\mu^{(2)}_{NQ_2}h & = \Lambda - T_3 \Lambda^3 + T_5 \Lambda^5 + \mathcal{O}(\Lambda)^7,
\end{aligned}
\end{equation}
where 
 \begin{equation} \label{eq:nq2igadet}
\begin{aligned}
T_3 & =  \frac{ 6 n_1^2 \varpi_1-6 n_1 \varpi_1+6 n_2^2 \varpi_2-6 n_2 \varpi_2+\varpi_2+\varpi_1 }{12 (\varpi_2+\varpi_1)}, \\
T_5  & =  \frac{1}{1440 (\varpi_2+\varpi_1)^2} \Big(  5 (6 n_1^2 \varpi_1-6 n_1 \varpi_1+6 n_2^2 \varpi_2-6 n_2 \varpi_2+\varpi_2+\varpi_1 )^2 \\ 
& \qquad \qquad \qquad \qquad \quad  + (\varpi_2+\varpi_1) \big( (180 n_1^4-360 n_1^3+120 n_1^2+60 n_1-17 ) \varpi_1 \\
& \qquad \qquad \qquad \qquad \quad  +(180 n_2^4-360 n_2^3+120 n_2^2+60 n_2-17 ) \varpi_2 \big) \Big).
\end{aligned}
\end{equation}

We seek a quadrature rule that reduces the dispersion error as much as possible. Since we consider a two-point rule, we have four degrees of freedom (two nodes and weights). To integrate exactly the $C_3^0$ space, the rule has to integrate the basis which, taking into account a repetitive pattern on uniform elements, involves three basis functions. This leads to the following algebraic system

\begin{equation} \label{eq:c2iga}
\begin{aligned}
T_3 & =  0, \\
T_5 & =  0, \\
3 n_1(1 - n_1)^2 \varpi_1 + 3 n_2 (1- n_2)^2 \varpi_2 & = \frac{1}{4}, \\
3 n_1^2(1 - n_1) \varpi_1 + 3 n_2^2 (1- n_2) \varpi_2 & = \frac{1}{4}, \\
n_1^3 \varpi_1 + n_2^3 \varpi_2 + (1- n_1)^3 \varpi_1 + (1 - n_2)^3 \varpi_2 & = \frac{1}{2}.
\end{aligned} 
\end{equation}

The system is built over the unit interval. The first two equations correspond to the dispersion error-minimization, while the last three represent the exactness of the rule on the $C_3^0$ spline space. 
There are five equations and four unknowns. 
However, symbolic calculations show that the equation $T_3 = 0$ is a redundant equation and this system has the following equivalent solutions:

Solution 1
\begin{equation} \label{eq:c2igasol1}
\begin{aligned}
n_1 & = \frac{1}{10} \Big(5 - \sqrt{ \frac{1}{3} (33 - 2 \sqrt{266} ) } \Big), \\
n_2 & = \frac{1}{150} \Big(75 - \sqrt{3}( 33 - 2 \sqrt{266} )^{3/2} + 66 \sqrt{ 3 (33 - 2 \sqrt{266} ) } \Big), \\
\varpi_1 & = \frac{1}{266} \big(133  +  2 \sqrt{266} \big), \\
\varpi_2 & = \frac{1}{266} \big(133  - 2 \sqrt{266} \big), \\
\end{aligned}
\end{equation}

Solution 2
\begin{equation} \label{eq:c2igasol2}
\begin{aligned}
n_1 & = \frac{1}{10} \Big(5 + \sqrt{ \frac{1}{3} (33 - 2 \sqrt{266} ) } \Big), \\
n_2 & = \frac{1}{150} \Big(75 + \sqrt{3}( 33 - 2 \sqrt{266} )^{3/2} - 66 \sqrt{ 3 (33 - 2 \sqrt{266} ) } \Big), \\
\varpi_1 & = \frac{1}{266} \big(133  +  2 \sqrt{266} \big), \\
\varpi_2 & = \frac{1}{266} \big(133  - 2 \sqrt{266} \big), \\
\end{aligned}
\end{equation}

Solution 3
\begin{equation} \label{eq:c2igasol3}
\begin{aligned}
n_1 & = \frac{1}{10} \Big(5 - \sqrt{ \frac{1}{3} (33 - 2 \sqrt{266} ) } \Big), \\
n_2 & = \frac{1}{150} \Big(75 - \sqrt{3}( 33 + 2 \sqrt{266} )^{3/2} + 66 \sqrt{ 3 (33 + 2 \sqrt{266} ) } \Big), \\
\varpi_1 & = \frac{1}{266} \big(133  -  2 \sqrt{266} \big), \\
\varpi_2 & = \frac{1}{266} \big(133  + 2 \sqrt{266} \big), \\
\end{aligned}
\end{equation}

Solution 4
\begin{equation} \label{eq:c2igasol4}
\begin{aligned}
n_1 & = \frac{1}{10} \Big(5 + \sqrt{ \frac{1}{3} (33 - 2 \sqrt{266} ) } \Big), \\
n_2 & = \frac{1}{150} \Big(75 + \sqrt{3}( 33 - 2 \sqrt{266} )^{3/2} - 66 \sqrt{ 3 (33 - 2 \sqrt{266} ) } \Big), \\
\varpi_1 & = \frac{1}{266} \big(133  -  2 \sqrt{266} \big), \\
\varpi_2 & = \frac{1}{266} \big(133  + 2 \sqrt{266} \big). \\
\end{aligned}
\end{equation}

These two-point quadrature rules lead to the discrete dispersion relation 
 \begin{equation} \label{eq:quadratico2dr}
 \begin{aligned} 
2( -\frac{1}{6} + \frac{7}{720} \Lambda^2 ) \cos(2\mu^{(2)}_{NQ_2}h) - 2( \frac{1}{3} - \frac{19}{90} \Lambda^2 ) \cos(\mu^{(2)}_{NQ_2}h) + (1 - \frac{67}{120} \Lambda^2 ) = 0,
 \end{aligned}
 \end{equation}
which gives the optimal dispersion error
\begin{equation} \label{eq:quadraticoptdre}
\mu^{(2)}_{NQ_2}h = \Lambda - \frac{11}{120960} \Lambda^7 - \frac{1}{345600} \Lambda^9 + \mathcal{O}(\Lambda)^{11}.
\end{equation}

\begin{remark} \label{rem:onep}
The last three equations in \eqref{eq:c2iga} define a one-parameter family of quadrature rules that integrates exactly the $C_3^0$ space. Among these rules, we seek those that further reduce the dispersion error. Such rules result in an approximation error of order seven, see \eqref{eq:quadraticoptdre}.
\end{remark}

\begin{remark} \label{rem:2iga}
Taking the difference between $\mu^{(2)}_{NQ_2}h$  and  $\Lambda$ gives an error representation of the dispersion error, which is of order seven. Symbolic calculations show that a three-point (or more points) rule does not increase the convergence order in the dispersion error. The extra degrees of freedom obtained by utilizing more quadrature points can be used to integrate a larger space than the polynomial space $C^0_3$. We present this alternative in the following subsection. The error with an order of seven as in \eqref{eq:quadraticoptdre} is the minimized dispersion error for $C^1$ quadratic isogeometric analysis.
\end{remark}

\begin{remark} \label{rem:2iga2}
This optimized dispersion error expansion coincides with the expansion derived from the blending schemes proposed in \cite{calo2017dispersion}. Further calculation shows that both the blending schemes and these new two-point quadrature rules give the same discrete dispersion relation.
\end{remark}

\subsection{A 2.5-point rule}
The two-point quadrature rules \eqref{eq:c2igasol1} to \eqref{eq:c2igasol4} integrate the space $C^0_3$ exactly. Alternatively, one can derive a 2.5-point rule which integrates the space $C_3^{-1}$ exactly as well as minimizes the dispersion errors. We denote $G_{2.5}$ a three-point quadrature rule with one of the points fixed at  one of the element boundaries, that might be shared with another element. That is, this class  of quadratures generalizes the Gauss-Radau family of quadratures. We denote this  quadrature rule $G_{2.5}$ for $C^1$ quadratic isogeometric elements with nodes $n_1, n_2, n_3$ in the reference interval $[0, 1]$ and weights $\varpi_1, \varpi_2, \varpi_3$. By setting $n_3=1$, we say this is a 2.5-point rule  per element as the mapping is continuous across the element interface, thus the evaluation at the interface does not need to be recomputed from the neighboring element. By the same process we described to solve \eqref{eq:c2iga}, we obtain the following quadrature rule.

\begin{equation} \label{eq:c2igasol5}
\begin{aligned}
n_1 & = \frac{1}{30} \Big(9 - \sqrt{ 51 } \Big), \\
n_2 & = \frac{1}{30} \Big(9 + \sqrt{ 51 } \Big), \\
n_3 & = 1, \\
\varpi_1 & = \frac{1}{442} \Big( 79  + 12 (9 - \sqrt{51} ) \Big), \\
\varpi_2 & = \frac{1}{442} \Big( 295  - 12 (9 - \sqrt{51} ) \Big), \\
\varpi_3 & = \frac{2}{13}.
\end{aligned}
\end{equation}

Alternatively, one can fix a quadrature point at $n_3 = 0$ and derive an alternative, but equivalent rule.

\begin{remark} \label{rem:2iga3}
For multidimensional case, we assume that a tensor product grid is placed on the domain $\Omega$. Then, we conclude that the above derivations are independent of the spatial dimension and the same rule remains valid for each dimension; more details are referred to \cite{calo2017dispersion}. Moreover, a duality between dispersion analysis and spectrum analysis in error expansion form was established in \cite{calo2017dispersion}. 
\end{remark}

\begin{figure}[ht]
\centering
\includegraphics[height=7.0cm]{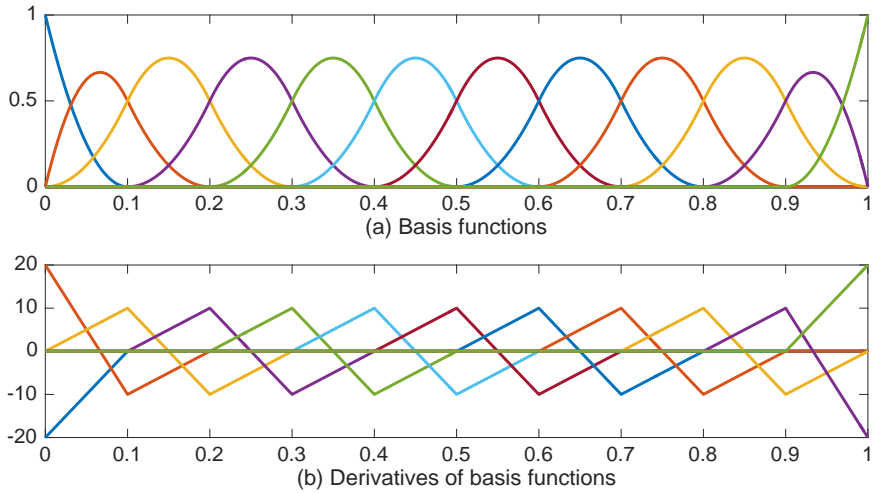}
\caption{Basis functions and their derivatives for $C^1$ quadratic isogeometric analysis.}
\label{fig:d}
\end{figure}

\subsection{Rule near the boundary elements}
Naturally, since the 2.5-point rule integrates the $C^{-1}_3$ space exactly, it can be applied to both the interior and boundary elements over the domain. This is different for the two-point rule as it takes advantage of the $C^0$ continuity of the integrand across element interfaces. This continuity assumption is lost at the boundary elements. Figure \ref{fig:d} (b) shows that the derivatives of the first (last) two basis function do not vanish at $x=0$ ($x=1$) associated with the left (right) boundary element, as is the case of all other element interfaces in the domain. 

Special treatment near these boundary elements is required. Simply, we propose to apply the 2.5-point rule at the boundary elements. This can be done in the implementation efficiently.

\section{Numerical experiments} \label{sec:num}

In this section, we present the numerical tests of the problem \eqref{eq:pde} in one and two dimensions (denoted with 1D, 2D) with uniform meshes. The comparisons of the numerical eigenvalues and eigenfunctions obtained from isogeometric analysis with those from finite elements are studied in \cite{cottrell2006isogeometric,hughes2008duality,hughes2014finite,puzyrev2017dispersion} and significant advantages of isogeometric elements over finite elements are shown. In this section, we show the numerical results of isogeometric elements obtained by the proposed new rules as well as the comparisons with those obtained by the optimal-blending rules.

The  exact eigenvalues and eigenfunctions of the one dimensional problem \eqref{eq:pde} are
\begin{equation}
\lambda_j = j^2 \pi^2, \quad \text{and} \quad u_j = \sqrt{2} \sin(j\pi x), \quad j = 1, 2, \cdots,
\end{equation}
respectively, while those of the two dimension problems are
\begin{equation}
\lambda_{jk} = ( j^2 + k^2 ) \pi^2, \quad \text{and} \quad u_{jk} = 2 \sin(j\pi x)\sin(k\pi y), \quad j,k = 1, 2, \cdots,
\end{equation}
respectively. We sort the approximate eigenvalues both in one and two dimension in the ascending order. In the following figures, we present the eigenvalue (EV) errors  as well as the eigenfunction (EF) errors in both $L^2$-norm and energy norm.

There are mainly three different optimally-blended rules proposed in \cite{calo2017dispersion,puzyrev2017dispersion}: three-point Gauss-Legendre rule with three-point Gauss-Lobatto rule; two-point Gauss-Legendre rule with three-point Gauss-Lobatto rule, and three-point Gauss-Legendre rule with two-point Gauss-Legendre rule. For the comparison with our new rules, we choose the last one as it  requires fewer evaluation points, and we denote this one as the blending rule for the tests. 

\begin{figure}[ht]
\centering
\includegraphics[height=6.0cm]{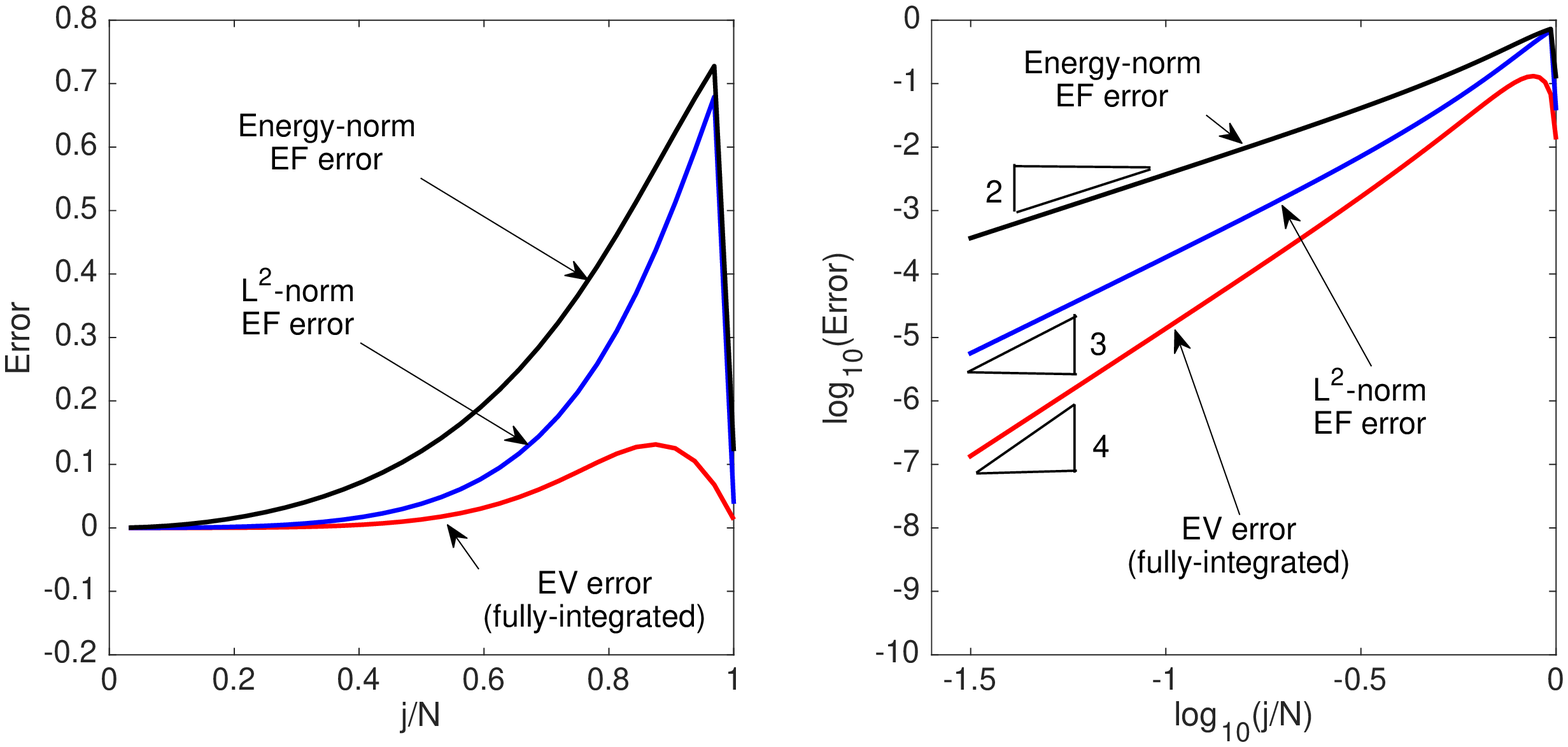} \\
 Fully-integrated inner products
\includegraphics[height=6.0cm]{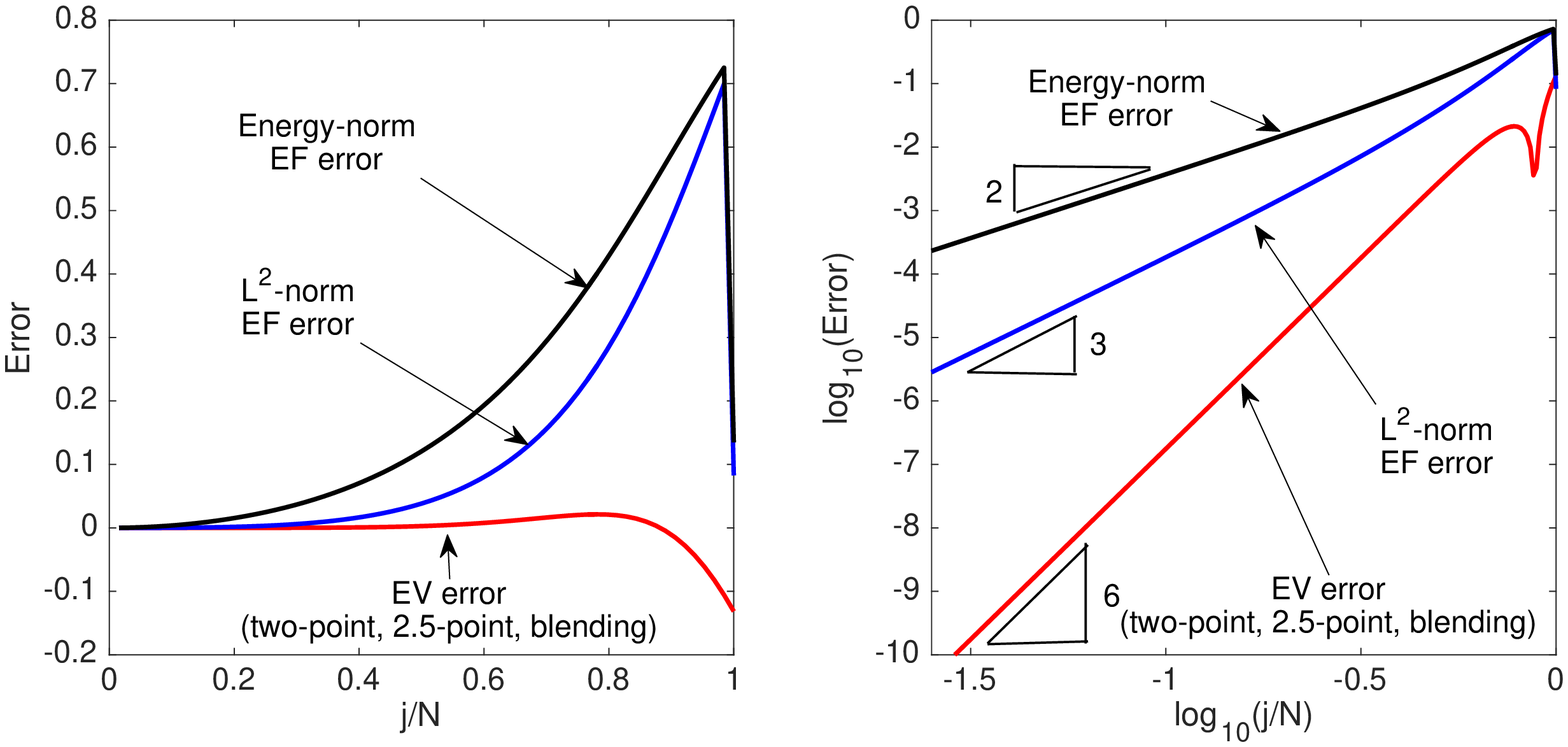} \\
Dispersion-optimized two-point, 2.5-point, and blending rules 
\caption{Eigenvalue (EV) and eigenfunction (EF) errors in linear (left) and logarithmic (right) scales using full integration (top) and dispersion optimized two-point, 2.5-point, and optimal blending rules (bottom) for $C^1$ quadratic isogeometric analysis.}
\label{fig:1d32grtpr}
\end{figure}

\begin{figure}[ht]
\centering
\includegraphics[height=6.2cm]{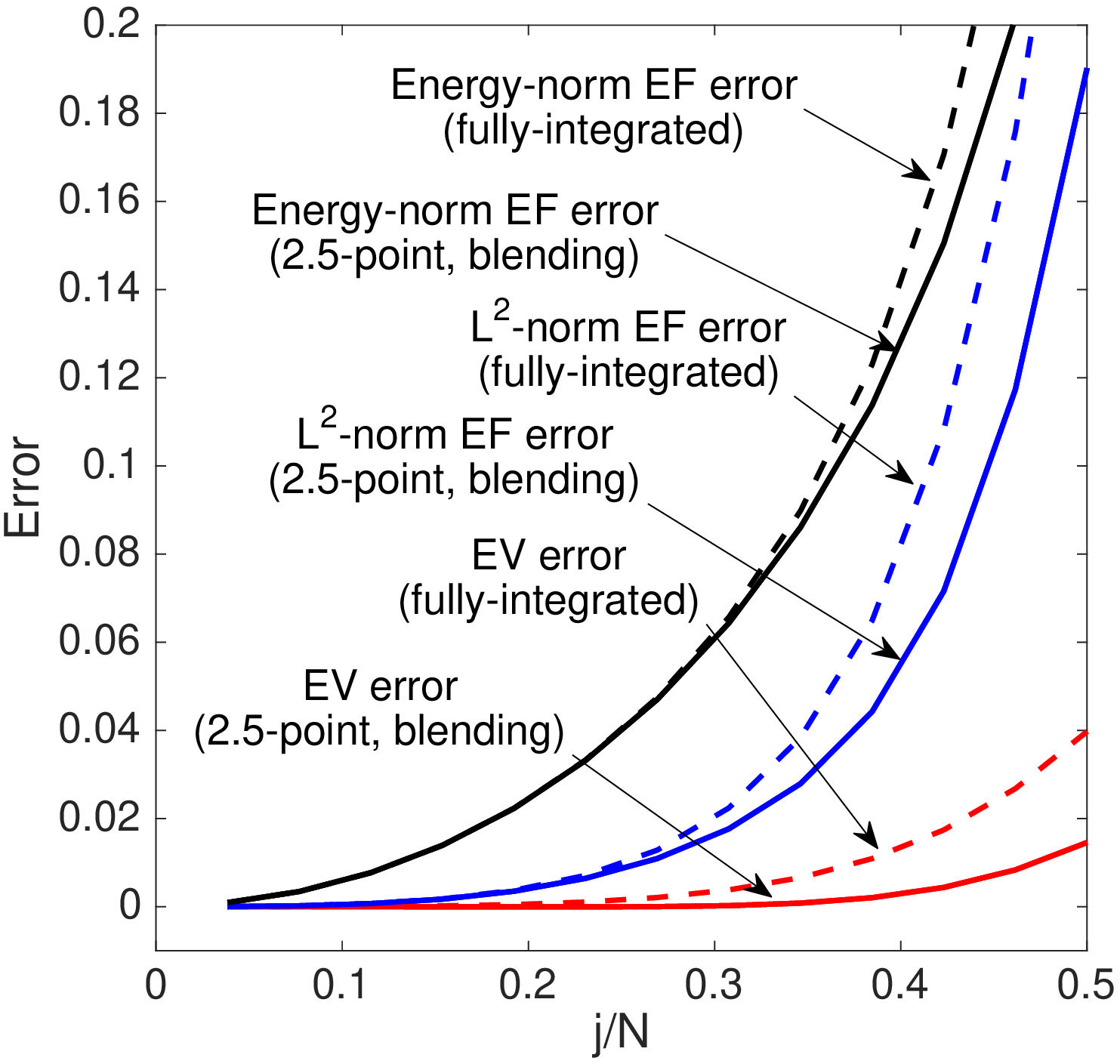}
\includegraphics[height=6.2cm]{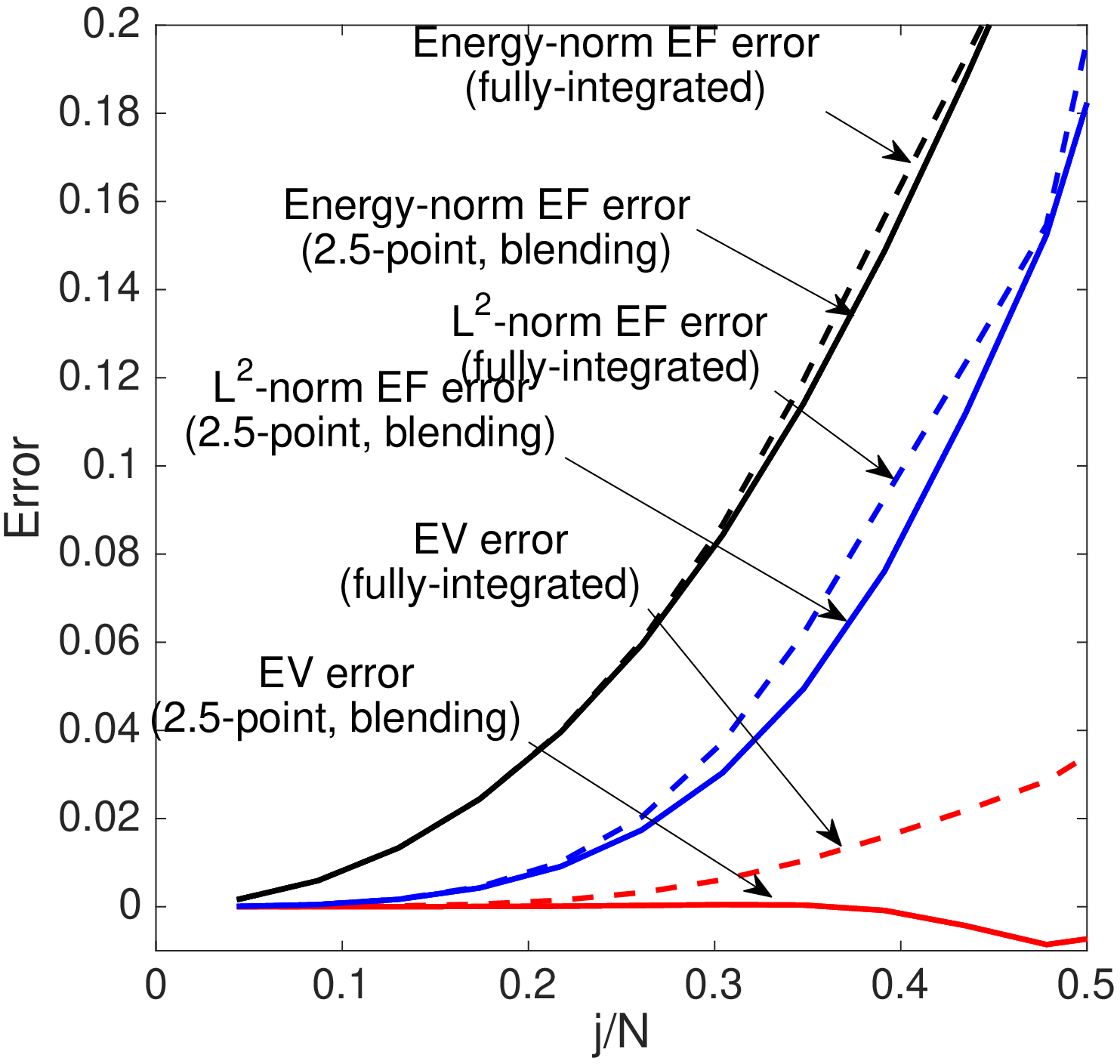}
\caption{Eigenvalue (EV) and eigenfunction (EF) errors using full integration, 2.5-point, and the optimal blending rules for $C^1$ quadratic isogeometric analysis on the stretched grid with stretching factor $1.05$ (left) and 1.07 (right).}
\label{fig:nu}
\end{figure}

In the one dimensional case, we consider the isogeometric elements with fully-integrated inner products and quadrature-rule approximated inner products. These quadrature rules include the two-point rule, 2.5-point rule, and the blending rule.

Figure \ref{fig:1d32grtpr} shows the relative eigenvalue errors $\dfrac{\lambda_j^h - \lambda_j}{\lambda_j}$, the $L^2$ eigenfunction errors $\| u_j^h - u_j \|_0$, and the scaled energy-norm errors $\dfrac{\| u_j^h - u_j \|_E}{\sqrt{\lambda_j}}$ in both linear and logarithmic scales when full integration and the  two-point, 2.5-point, and the blending rules are applied for the inner product for $C^1$ quadratic isogeometric elements. We scale the energy-norm error in the view of the generalized Pythagorean eigenvalue error theorem \cite{puzyrev2017dispersion,calo2017dispersion}. In the implementation, for this example, the full integration is realized by applying the three-point Gauss rule.  In Figure \ref{fig:1d32grtpr}, we observe that two extra orders of convergence in the eigenvalue errors when the two-point, 2.5-point, and the blending rules are applied. Also, we observe that the two-point rule, 2.5-point rule, and the optimal-blending rule lead to the same results. This verifies numerically that these newly-developed rules are equivalent to the optimally-blended rules.

\begin{figure}[ht]
\centering
\includegraphics[height=7.2cm]{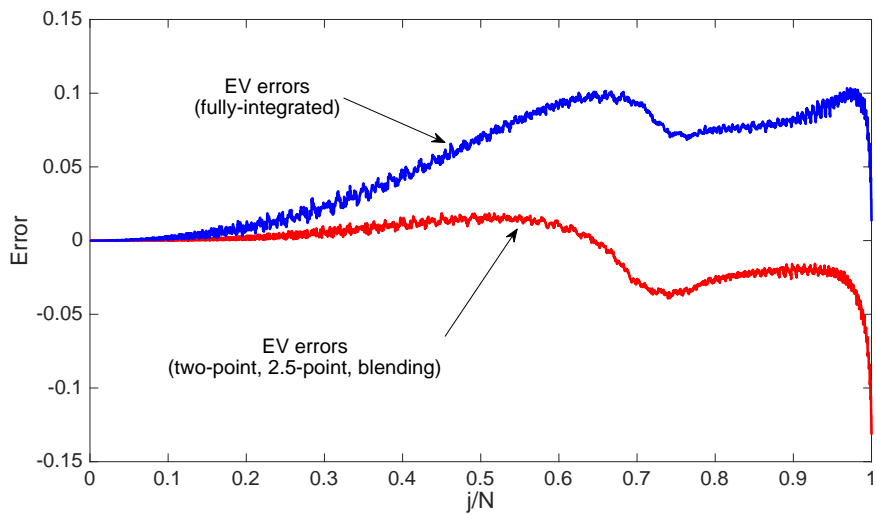}
\caption{Eigenvalue (EV) errors in 2D using full integration, two-point, 2.5-point, and optimal blending rules for $C^1$ quadratic isogeometric analysis.}
\label{fig:2dcomerr}
\end{figure}

A similar behavior is observed on non-uniform meshes as shown in Figure \ref{fig:nu} and in 2D as shown in Figure \ref{fig:2dcomerr}.
Figure \ref{fig:nu} shows an equivalence between the 2.5-point and the blending rules for $C^1$ quadratic isogeometric analysis on non-uniform meshes, precisely the stretching meshes with stretching factors 1.05 and 1.07. 
A simple 2D test example on a uniform $64 \times 64$ mesh is shown in Figure \ref{fig:2dcomerr}. We observe that both two-point and 2.5-point rule lead to the same results as those from optimally-blended rules. For more numerical results in 2D, we refer to the paper \cite{puzyrev2017dispersion}.

\section{Concluding remarks} \label{sec:conclusion}
We introduce new quadrature rules that minimize the dispersion error and possess the optimal septic convergence order. The optimal rule requires two nodes per element and arises from an algebraic system that admits a closed-form solution. This rule is valid for uniform grids with periodic boundary conditions. We also introduce a 2.5-rule which exactly integrates discontinuous, cubic functions on general grids. The optimal rule combined with the generalized 2.5-point Gauss-Radau type of rule on the boundary elements remains dispersion-minimizing on finite domains with arbitrary boundary conditions. Moreover, compared with the optimal blending schemes proposed in \cite{calo2017dispersion}, our approach further reduces the number of quadrature points, which brings a significant computational speed-up to the application problems such as the wave propagation or structural vibrations, particularly in three dimensions.  Future work in this direction includes further studies on non-uniform meshes as well as extension to higher order isogeometric analysis.

\section{Acknowledgments} 
This publication was made possible by the CSIRO Professorial Chair in Computational Geoscience of Curtin University, additional support was provided by  the National Priorities Research Program grant 7-1482-1-278 from the Qatar National Research Fund (a member of The Qatar Foundation), and by the European Union's Horizon 2020 Research and Innovation Program of the Marie Sk{\l}odowska-Curie grant agreement No. 644202. The Spring 2016 Trimester on "Numerical methods for PDEs", organized  with the collaboration of the Centre Emile Borel at the Institut Henri Poincare in Paris supported VMC's visit to IHP in October, 2016.

\section*{References}

\bibliographystyle{siam}


\bibliography{igaref}

\end{document}